\documentclass[12pt]{article}
\usepackage{mathrsfs}
\usepackage{amssymb} \textwidth 150mm \textheight 220mm \oddsidemargin
15pt \evensidemargin 0pt \topmargin 0cm \headsep 0.3cm

\usepackage{amsmath}
\usepackage{graphicx}
\usepackage{indentfirst}
\usepackage{cite}
\usepackage{float}
\usepackage{CJK}

\begin{document}
\title{\Large\bf{Bifurcation of limit cycles from a quadratic global center with two switching lines
\thanks{E-mail address: jihua1113@163.com (J. Yang)} }}
\author{{Jihua Yang}\\
{\small \it School of Mathematics and Computer Science, Ningxia Normal University,}\\
 {\small \it  Guyuan, 756000, China}
 }
\date{}
\maketitle \baselineskip=0.9\normalbaselineskip \vspace{-3pt}
\noindent
{\bf Abstract}\, In this paper, we generalize the Picard-Fuchs equation method to study the bifurcation of limit cycles of perturbed differential systems with two switching lines. We obtain the detailed expression of the corresponding first order Melnikov function which can be used to get the upper bound of the number of limit cycles for the perturbed system by using Picard-Fuchs equation. It is worth noting that we greatly simplify the computations and this method can be applied to study the number of limit cycles of other differential systems with two switching lines. Our results also show that the number of switching lines has essentially impact on the number of limit cycles bifurcating from a period annulus.
\vskip 0.2 true cm
\noindent

 \section{Introduction and main results}
 \setcounter{equation}{0}
\renewcommand\theequation{1.\arabic{equation}}

Piecewise differential systems belong to an interesting class of nonlinear differential systems to be studied  extensively as they frequently appear in modeling many real phenomena. For instance, in control engineering \cite{BBCK}, nonlinear oscillations \cite{T} and biology \cite{K}, etc.. Moreover, these systems can exhibit complicated dynamical phenomena such as those exhibited by general nonlinear differential systems. Thus, these past years, much interest from the mathematical community is seen in trying to understand their dynamical richness, especially the bifurcation of limit cycles.

There are many excellent papers studying the bifurcation of limit cycles of piecewise differential systems with one switching line, see for example, \cite{LH,ZY,YZ17,YZ18,LHR,LM,LMN,DL,CLZ,LL,GPL} and the references quoted there. Of course, there are also a few papers dedicated to study bifurcation of limit cycles of piecewise differential systems with multiple switching lines, see \cite{IL,CLXY,BLV,HD,SD,LYD}. The methods used in the above papers are main the Melnikov function established in \cite{LH,HS} and  the averaging method developed in \cite{BL,LMN,SV}. The disadvantage of the above two methods is the complexity of involved calculations. Recently, Yang and Zhao \cite{YZ18} developed the Picard-Fuchs equation method to study the number of limit cycles of piecewise differential systems with one switching line.

In this paper, our aim is to generalize the Picard-Fuchs equation method to study the bifurcation of limit cycles of perturbed differential systems with two switching lines. More precisely, we study the following integrable differential system under the perturbations of piecewise polynomials of degree $n$
\begin{eqnarray}
       \dot{x}=y-2x^2-\eta, \ \ \
        \dot{y}=-2xy,
\end{eqnarray}
where $\eta$ is a real positive constant. System (1.1)  has a unique global center $G(0,\eta)$. 
The perturbed system of (1.1) with two vertical switching lines intersected at point $(0,\eta)$ is
\begin{eqnarray}
\left(  \begin{array}{c}
          \dot{x} \\          \dot{y}
          \end{array} \right)
=\begin{cases}
 \left(
  \begin{array}{c}
          y-2x^2-\eta+\varepsilon f^1(x,y) \\
          -2xy+\varepsilon g^1(x,y)
          \end{array} \right), \quad x>0,\ y>\eta,\\[0.6truecm]
 \left(  \begin{array}{c}
         y-2x^2-\eta+\varepsilon f^2(x,y) \\
          -2xy+\varepsilon g^2(x,y)
           \end{array}
 \right),\quad x>0,\ y<\eta,\\[0.6truecm]
 \left(  \begin{array}{c}
         y-2x^2-\eta+\varepsilon f^3(x,y) \\
          -2xy+\varepsilon g^3(x,y)
           \end{array}
 \right),\quad x<0,\ y<\eta,\\[0.6truecm]
 \left(  \begin{array}{c}
         y-2x^2-\eta+\varepsilon f^4(x,y) \\
          -2xy+\varepsilon g^4(x,y)
           \end{array}
 \right),\quad x<0,\ y>\eta,
  \end{cases}
 \end{eqnarray}
 where $0<|\varepsilon|\ll1$,
$$f^k(x,y)=\sum\limits_{i+j=0}^na^k_{i,j}x^iy^j,\ \ g^k(x,y)=\sum\limits_{i+j=0}^nb^k_{i,j}x^iy^j,\ k=1,2,3,4.$$
When $\varepsilon=0$, system (1.2) has first integrals
\begin{eqnarray}
\begin{aligned}
&H^1(x,y)=y^{-2}\Big(x^2-y+\frac{\eta}{2}\Big)=h, \ x>0,\ y>\eta,\\
&H^2(x,y)=y^{-2}\Big(x^2-y+\frac{\eta}{2}\Big)=h, \ x>0,\ y<\eta,\\
&H^3(x,y)=y^{-2}\Big(x^2-y+\frac{\eta}{2}\Big)=h, \ x<0,\ y<\eta,\\
&H^4(x,y)=y^{-2}\Big(x^2-y+\frac{\eta}{2}\Big)=h, \ x<0,\ y>\eta
\end{aligned}
\end{eqnarray}
with integrating factor $\mu^k(x,y)=y^{-3},\ k=1,2,3,4$ and a family of periodic orbits given by
\begin{eqnarray}
\begin{aligned}
L_h=&\{H^1(x,y)=h, x>0,\ y>\eta\}\cup\{H^2(x,y)=h, x>0,\ y<\eta\}\\
&\cup\{H^3(x,y)=h, x<0,\ y<\eta\}\cup\{H^4(x,y)=h, x<0,\ y>\eta\}\\
&:=L^1_h\cup L^2_h\cup L^3_h\cup L^4_h, \ h\in\Sigma=\Big(-\frac{1}{2\eta},0\Big).
\end{aligned}
\end{eqnarray}
 Obviously, $L_h$ approaches the center $G(0,\eta)$ as $h\rightarrow-\frac{1}{2\eta}$ and an invariant curve $x^2-y+\frac{\eta}{2}=0$ as $h\rightarrow0$, respectively. See Fig.\,1.

 Our main results are as follows.
\vskip 0.2 true cm

\noindent
 {\bf Theorem 1.1.} {\it The number of limit cycles of system (1.2) bifurcating from the period annulus around the global center is not more than $41n-23$ (counting multiplicity) for $n=1,2,3,\cdots$.}
 \vskip 0.2 true cm

\noindent
 {\bf Theorem 1.2.} {\it If $f^1(x,y)=f^2(x,y)$ and $f^3(x,y)=f^4(x,y)$, then the number of limit cycles of system (1.2) bifurcating from the period annulus around the global center is not more than $9n-4$ (counting multiplicity) for $n=1,2,3,\cdots$.}
  \vskip 0.2 true cm

\noindent
 {\bf Theorem 1.3.} {\it If $f^1(x,y)=f^4(x,y)$ and $f^2(x,y)=f^3(x,y)$, then the number of limit cycles of system (1.2) bifurcating from the period annulus around the global center is not more than $9n-6$ (counting multiplicity) for $n=1,2,3,\cdots$.}
 \vskip 0.2 true cm

\noindent
{\bf Remark 1.1}. When $f^1(x,y)=f^2(x,y)=f^3(x,y)=f^4(x,y)$ and $g^1{x,y}=g^2(x,y)=g^3(x,y)=g^4(x,y)$, Gentes \cite{G} studied the case for $n=2$ and proved that $M(h)$ has at most 2 zeros. Xiong and Han \cite{XH} studied the case for $n$ and obtained that $M(h)$ has at most $n$ zeros.
\vskip 0.2 true cm

\noindent
{\bf Remark 1.2}. From Theorems 1.1-1.3, we know that the number of switching lines has essentially impact on the number of limit cycles bifurcating from the global center.
\vskip 0.2 true cm

\noindent
{\bf Remark 1.3}. The techniques we use mainly include Melnikov function, Picard-Fuchs equation and derivation-division algorithm. We first obtain the generators of the first order Melnikov function $M(h)$ and the Picard-Fuchs equations which they satisfy. Thus, we can compute the generators by using the Picard-Fuchs equations. Finally, we get the number of the zeros of $M(h)$ by derivation-division algorithm. It is worth noting that we greatly simplify the computations by Picard-Fuchs equations and this method can be applied to study the number of limit cycles of other differential systems with two switching lines.
\vskip 0.2 true cm

The rest of the paper is organized as follows: In Section 2, we will give detailed expression of the first Melnikov function $M(h)$ by using Picard-Fuchs equation. Theorem 1.1-1.3 will be proved in Sections 3-5.

\section{The algebraic structure of $M(h)$ and Picard-Fuchs equation}
 \setcounter{equation}{0}
\renewcommand\theequation{2.\arabic{equation}}

Noting that the direction of the periodic orbits is clockwise, by the Theorem 2.2 in \cite{HD} and the Lemma 2.1 in \cite{WHC}, we know that the first order Melnikov function $M(h)$ of system (1.2) has the following form
\begin{eqnarray*}
\begin{aligned}
M(h)=&\frac{H^1_y(A)H^2_x(B)H^3_y(C)H^4_x(D)}{H^4_y(A)H^1_x(B)H^2_y(C)H^3_x(D)}
\int_{L^1_h}\mu^1(x,y)[g^1(x,y)dx-f^1(x,y)dy]\\
&+\frac{H^1_y(A)H^3_y(C)H^4_x(D)}{H^4_y(A)H^2_y(C)H^3_x(D)}
\int_{L^2_h}\mu^2(x,y)[g^2(x,y)dx-f^2(x,y)dy]\\
&+\frac{H^1_y(A)H^4_x(D)}{H^4_y(A)H^3_x(D)}
\int_{L^3_h}\mu^3(x,y)[g^3(x,y)dx-f^3(x,y)dy]\\
&+\frac{H^1_y(A)}{H^4_y(A)}
\int_{L^4_h}\mu^4(x,y)[g^4(x,y)dx-f^4(x,y)dy],\ h\in\Sigma
\end{aligned}
\end{eqnarray*}
and the number of zeros of $M(h)$ controls the number of limit cycles of system (1.2) if $M(h)\not\equiv0$
in the corresponding period annulus. It is easy to check that
$$\frac{H^1_y(A)H^2_x(B)H^3_y(C)H^4_x(D)}{H^4_y(A)H^1_x(B)H^2_y(C)H^3_x(D)}
=\frac{H^1_y(A)H^3_y(C)H^4_x(D)}{H^4_y(A)H^2_y(C)H^3_x(D)}=
\frac{H^1_y(A)H^4_x(D)}{H^4_y(A)H^3_x(D)}=\frac{H^1_y(A)}{H^4_y(A)}=1.$$
Hence,
\begin{eqnarray}
\begin{aligned}
M(h)=&\int_{L^1_h}y^{-3}[g^1(x,y)dx-f^1(x,y)dy]+\int_{L^2_h}y^{-3}[g^2(x,y)dx-f^2(x,y)dy]\\
&+\int_{L^3_h}y^{-3}[g^3(x,y)dx-f^3(x,y)dy]+\int_{L^4_h}y^{-3}[g^4(x,y)dx-f^4(x,y)dy].
\end{aligned}
\end{eqnarray}

For $h\in\Sigma$ and $i=0,1,2,\cdots,j=-1,0,1,2,\cdots$, we denote
$$\begin{aligned}
&I_{i,j}(h)=\int_{L^1_h}x^iy^{j-3}dy,\ \ J_{i,j}(h)=\int_{L^2_h}x^iy^{j-3}dy,\\
&\tilde{J}_{i,j}(h)=\int_{L^3_h}x^iy^{j-3}dy,\ \ \tilde{I}_{i,j}(h)=\int_{L^4_h}x^iy^{j-3}dy.
\end{aligned}$$
Let $\Omega$ be the interior of $L^1_{h}\cup \overrightarrow{BG}\cup \overrightarrow{GA}$, see the black line in Fig.\,1. Using the Green's Formula, we have for $i\geq0$ and $j\geq-1$
\begin{eqnarray*}
\begin{aligned}
\int_{L^1_{h}}x^iy^jdx
&=\oint_{L^1_{h}\cup\overrightarrow{BG}\cup \overrightarrow{GA}}x^iy^jdx-\int_{\overrightarrow{BG}}x^iy^jdx\\
&=j\iint\limits_\Omega x^{i}y^{j-1}dxdy-\eta^j\int_{\overrightarrow{BG}}x^idx,\\
\int_{L^1_{h}}x^{i+1}y^{j-1}dy&=\oint_{L^1_{h}\cup\overrightarrow{BG}\cup \overrightarrow{GA}}x^{i+1}y^{j-1}dy
=-(i+1)\iint\limits_\Omega x^{i}y^{j-1}dxdy.
\end{aligned}
\end{eqnarray*}
Hence,
\begin{eqnarray}
\begin{aligned}
\int_{L^1_{h}}x^iy^jdx=-\frac{j}{i+1}\int_{L^1_h}x^{i+1}y^{j-1}dy-\eta^j\int_{\overrightarrow{BG}}x^idx.
\end{aligned}
\end{eqnarray}
In a similar way, we have for $i\geq 0$ and $j\geq-1$
\begin{eqnarray}
\begin{aligned}
\int_{L^2_{h}}x^iy^jdx=-\frac{j}{i+1}\int_{L^2_h}x^{i+1}y^{j-1}dy-\eta^j\int_{\overrightarrow{GB}}x^idx,\\
\int_{L^3_{h}}x^iy^jdx=-\frac{j}{i+1}\int_{L^3_h}x^{i+1}y^{j-1}dy-\eta^j\int_{\overrightarrow{DG}}x^idx,\\
\int_{L^4_{h}}x^iy^jdx=-\frac{j}{i+1}\int_{L^4_h}x^{i+1}y^{j-1}dy-\eta^j\int_{\overrightarrow{GD}}x^idx.
\end{aligned}
\end{eqnarray}
Therefore, we have from (2.1)-(2.3)
\begin{eqnarray*}
\begin{aligned}
M(h)=&\sum\limits_{i+j=0}^nb^1_{i,j}\int_{L^1_h}x^iy^{j-3}dx
-\sum\limits_{i+j=0}^na^1_{i,j}\int_{L_h^1}x^iy^{j-3}dy\\
&+\sum\limits_{i+j=0}^nb^2_{i,j}\int_{L^2_h}x^iy^{j-3}dx
-\sum\limits_{i+j=0}^na^2_{i,j}\int_{L_h^2}x^iy^{j-3}dy\\
&+\sum\limits_{i+j=0}^nb^3_{i,j}\int_{L^3_h}x^iy^{j-3}dx
-\sum\limits_{i+j=0}^na^3_{i,j}\int_{L_h^3}x^iy^{j-3}dy\\
&+\sum\limits_{i+j=0}^nb^4_{i,j}\int_{L^4_h}x^iy^{j-3}dx
-\sum\limits_{i+j=0}^na^4_{i,j}\int_{L_h^4}x^iy^{j-3}dy\\
=&-\sum\limits_{i+j=0}^nb^1_{i,j}\Big(\frac{j-3}{i+1}\int_{L_h^1}x^{i+1}y^{j-4}dy+\eta^{j-3}\int_{\overrightarrow{BG}}x^idx\Big)
-\sum\limits_{i+j=0}^na^1_{i,j}\int_{L^1_h}x^iy^{j-3}dy\\
&-\sum\limits_{i+j=0}^nb^2_{i,j}\Big(\frac{j-3}{i+1}\int_{L_h^2}x^{i+1}y^{j-4}dy+\eta^{j-3}\int_{\overrightarrow{GB}}x^idx\Big)
-\sum\limits_{i+j=0}^na^2_{i,j}\int_{L^2_h}x^iy^{j-3}dy\\
&-\sum\limits_{i+j=0}^nb^3_{i,j}\Big(\frac{j-3}{i+1}\int_{L_h^3}x^{i+1}y^{j-4}dy+\eta^{j-3}\int_{\overrightarrow{DG}}x^idx\Big)
-\sum\limits_{i+j=0}^na^3_{i,j}\int_{L^3_h}x^iy^{j-3}dy\\
&-\sum\limits_{i+j=0}^nb^4_{i,j}\Big(\frac{j-3}{i+1}\int_{L_h^4}x^{i+1}y^{j-4}dy+\eta^{j-3}\int_{\overrightarrow{GD}}x^idx\Big)
-\sum\limits_{i+j=0}^na^4_{i,j}\int_{L^4_h}x^iy^{j-3}dy\\
=&\sum\limits_{\substack{i+j=0,\\i\geq0,j\geq-1}}^n\tilde{a}_{i,j}I_{i,j}(h)
+\sum\limits_{i=0}^n\tilde{a}_{i}\int_{\overrightarrow{BG}}x^idx
+\sum\limits_{\substack{i+j=0,\\i\geq0,j\geq-1}}^n\tilde{b}_{i,j}J_{i,j}(h)
+\sum\limits_{i=0}^n\tilde{b}_{i}\int_{\overrightarrow{GB}}x^idx\\
&+\sum\limits_{\substack{i+j=0,\\i\geq0,j\geq-1}}^n\tilde{c}_{i,j}\tilde{J}_{i,j}(h)
+\sum\limits_{i=0}^n\tilde{c}_{i}\int_{\overrightarrow{DG}}x^idx
+\sum\limits_{\substack{i+j=0,\\i\geq0,j\geq-1}}^n\tilde{d}_{i,j}\tilde{I}_{i,j}(h)
+\sum\limits_{i=0}^n\tilde{d}_{i}\int_{\overrightarrow{GD}}x^idx,\\
=&\sum\limits_{\substack{i+j=0,\\i\geq0,j\geq-1}}^n\sigma_{i,j}I_{i,j}(h)+\sum\limits_{\substack{i+j=0,\\i\geq0,j\geq-1}}^n\tau_{i,j}J_{i,j}(h)\\
&+\sum\limits_{i=0}^n\tilde{a}_{i}\int_{\overrightarrow{BG}}x^idx
+\sum\limits_{i=0}^n\tilde{b}_{i}\int_{\overrightarrow{GB}}x^idx
+\sum\limits_{i=0}^n\tilde{c}_{i}\int_{\overrightarrow{DG}}x^idx
+\sum\limits_{i=0}^n\tilde{d}_{i}\int_{\overrightarrow{GD}}x^idx,
\end{aligned}
\end{eqnarray*}
where $\tilde{a}_{i,j}$, $\tilde{b}_{i,j}$, $\tilde{c}_{i,j}$, $\tilde{d}_{i,j}$, $\sigma_{i,j}$, $\tau_{i,j}$, $\tilde{a}_{i}$, $\tilde{b}_{i}$, $\tilde{c}_{i}$ and $\tilde{d}_{i}$ are arbitrary real constants and in the last equality we have used
$$\tilde{I}_{i,j}(h)=(-1)^{i+1}I_{i,j}(h),\  \tilde{J}_{i,j}(h)=(-1)^{i+1}J_{i,j}(h).$$
The coordinates of $B$ and $D$ are $(\sqrt{\eta^2h+\frac{\eta}{2}},\eta)$ and $(-\sqrt{\eta^2h+\frac{\eta}{2}},\eta)$ respectively.
Thus,
\begin{eqnarray}
\begin{aligned}
M(h)=\sum\limits_{\substack{i+j=0,\\i\geq0,j\geq-1}}^n\sigma_{i,j}I_{i,j}(h)+\sum\limits_{\substack{i+j=0,\\i\geq0,j\geq-1}}^n\tau_{i,j}J_{i,j}(h)+
\sum\limits_{i=0}^n\nu_i\eta^{i+1}\Big(h+\frac{1}{2\eta}\Big)^{\frac{i+1}{2}},
\end{aligned}
\end{eqnarray}
where $\nu_i$ is a real constant.
\vskip 0.2 true cm

\noindent
{\bf Lemma 2.1.}\, {\it If $i+j=n\geq3$ and $i$ is an even number, then
\begin{eqnarray}
\begin{aligned}
&I_{i,j}(h)=\frac{1}{h^{n-2}}\Big[\tilde{\alpha}_1(h)I_{0,1}(h)+\tilde{\beta}_1(h)I_{2,0}(h)+\tilde{\varphi}_{\frac{3}{2}n-\frac{7+(-1)^n}{4}}(h)\Big],\\
&J_{i,j}(h)=\frac{1}{h^{n-2}}\Big[\tilde{\alpha}_2(h)J_{0,1}(h)+\tilde{\beta}_2(h)J_{2,0}(h)+\tilde{\psi}_{\frac{3}{2}n-\frac{7+(-1)^n}{4}}(h)\Big].
\end{aligned}
\end{eqnarray}
If $i+j=n\geq3$ and $i$ is an odd number, then
\begin{eqnarray*}
\begin{aligned}
&I_{i,j}(h)=\frac{1}{h^{n-2}}\Big[\tilde{\gamma}_1(h)I_{1,0}(h)+\tilde{\delta}_1(h)I_{1,1}(h)+\sqrt{ h+\frac{1}{2\eta}}\bar{\varphi}_{\frac{3}{2}n-\frac{9-(-1)^n}{4}}(h)\Big],\\
&J_{i,j}(h)=\frac{1}{h^{n-2}}\Big[\tilde{\gamma}_2(h)J_{1,0}(h)+\tilde{\delta}_2(h)J_{1,1}(h)+\sqrt{ h+\frac{1}{2\eta}}\bar{\psi}_{\frac{3}{2}n-\frac{9-(-1)^n}{4}}(h)\Big],
\end{aligned}
\end{eqnarray*}
where $\tilde{\varphi}_{l}(h)$, $\tilde{\psi}_{l}(h)$, $\bar{\varphi}_{l}(h)$ and $\bar{\psi}_{l}(h)$ are polynomials in $h$ of degrees at most $l$, and $\tilde{\alpha}_k(h)$, $\tilde{\beta}_k(h)$, $\tilde{\gamma}_k(h)$ and $\tilde{\delta}_k(h)$ are polynomials of $h$  with
\begin{eqnarray*}
\begin{aligned}
&\deg\tilde{\alpha}_k(h)\leq n-\frac{3+(-1)^n}{2},\ \deg\tilde{\delta}_k(h)\leq n-\frac{3-(-1)^n}{2},\\
&\deg\tilde{\beta}_k(h),\deg\tilde{\gamma}_k(h)\leq n-2,\ k=1,2.
\end{aligned}\end{eqnarray*}}

 \vskip 0.2 true cm

\noindent
{\bf Proof.}\, Without loss of generality, we only prove the first equality in (2.5). The others can be shown in a similar way. It follows from (1.3) that
\begin{eqnarray}
-2x^2y^{-3}+2xy^{-2}\frac{\partial x}{\partial y}+y^{-2}-\eta y^{-3}=0.
\end{eqnarray}
Multiplying (2.6) by $x^{i-2}y^{j}dy$, integrating over $L^1_h$ and noting that (2.2), we have
\begin{eqnarray}
\begin{aligned}
2(i+j-2)I_{i,j}(h)=iI_{i-2,j+1}(h)-\eta iI_{i-2,j}(h)+2\eta^{i+j-2}\Big(h+\frac{1}{2\eta}\Big)^{\frac{i}{2}}.
\end{aligned}
\end{eqnarray}
Similarly, multiplying the first equality in (1.3) by $x^{i}y^{j-3}dx$ and integrating over $L^1_h$ yields
\begin{eqnarray}
hI_{i,j}(h)=I_{i+2,j-2}(h)-I_{i,j-1}(h)+\frac{\eta}{2}I_{i,j-2}(h).
\end{eqnarray}
Taking $(i,j)=(2,0),(3,-1)$ in (2.7), we obtain
\begin{eqnarray}
\begin{aligned}
&I_{0,0}(h)=\eta^{-1}I_{0,1}(h)+\eta^{-1}\Big( h+\frac{1}{2\eta}\Big),\\ &I_{1,-1}(h)=\eta^{-1}I_{1,0}(h)+\frac{2}{3}\eta^{-1}\Big( h+\frac{1}{2\eta}\Big)^\frac{3}{2}.
\end{aligned}
\end{eqnarray}
From (2.8) we obtain
\begin{eqnarray}
\begin{aligned}
&I_{0,2}(h)=\frac{1}{h}\Big(I_{2,0}(h)-I_{0,1}(h)+\frac{\eta}{2}I_{0,0}(h)\Big),\\
&I_{3,-1}(h)=hI_{1,1}(h)+I_{1,0}(h)-\frac{\eta}{2}I_{1,-1}(h).
\end{aligned}
\end{eqnarray}
Taking $(i,j)=(2,-1)$ in (2.7) and $(i,j)=(0,1)$ in (2.8), we have
\begin{eqnarray}
\begin{aligned}
&I_{2,-1}(h)=\eta I_{0,-1}(h)-I_{0,0}(h)-\eta^{-2}\Big(\eta h+\frac{1}{2}\Big),\\
&hI_{0,1}(h)=I_{2,-1}(h)-I_{0,0}(h)+\frac{\eta}{2}I_{0,-1}(h).
\end{aligned}
\end{eqnarray}
Eliminating $I_{0,-1}(h)$ in (2.11) and noting that (2.9), we get
\begin{eqnarray}
I_{2,-1}(h)=\frac{1}{3}(2h+\eta^{-1})I_{0,1}(h).
\end{eqnarray}
From (2.7) and (2.8) we have
\begin{eqnarray}
\left\{\begin{array}{l}
I_{0,3}(h)=\frac{1}{h}\Big(I_{2,1}(h)-I_{0,2}(h)+\frac{\eta}{2}I_{0,1}(h)\Big),\\
I_{1,2}(h)=\frac{1}{h}\Big(I_{3,0}(h)-I_{1,1}(h)+\frac{\eta}{2}I_{1,0}(h)\Big),\\
I_{2,1}(h)=I_{0,2}(h)-\eta I_{0,1}(h)+\eta h+\frac{1}{2},\\
I_{3,0}(h)=\frac{3}{2}I_{1,1}(h)-\frac{3}{2}\eta I_{1,0}(h)+\eta\Big(h+\frac{1}{2\eta}\Big)^\frac{3}{2},\\
I_{4,-1}(h)=2I_{2,0}(h)-2\eta I_{2,-1}(h)+h+\frac{1}{2\eta}.
\end{array}\right.
\end{eqnarray}

Now we prove the conclusion by induction on $n$. In fact, (2.13) implies that the conclusion holds for $n=3$. Suppose that the first equality in (2.5) holds for $i+j\leq n-1\, (n\geq4)$. If $n$ is an even number, then, by (2.7) and (2.8), we have
\begin{eqnarray}
\mathbf{A}\left(\begin{matrix}
                I_{0,n}(h)\\
                 I_{2,n-2}(h)\\
                 I_{4,n-4}(h)\\
                 \vdots\\
                  I_{n-2,2}(h)\\
                  I_{n,0}(h)
                \end{matrix}\right)\ \
=\left(\begin{matrix}
                \frac{1}{h}\big[-I_{0,n-1}(h)+\frac{\eta}{2}I_{0,n-2}(h)\big]\\
                \frac{1}{n-2}\big[I_{0,n-1}(h)-\eta I_{0,n-2}(h)+\eta^{n-2}(h+\frac{1}{2\eta})\big]\\
                \frac{1}{n-2}\big[2I_{2,n-3}(h)-2\eta I_{2,n-4}(h)+\eta^{n-2}(h+\frac{1}{2\eta})^2\big]\\
                                  \vdots\\
                 \frac{1}{2n-4}\big[(n-2)I_{n-4,3}(h)-(n-2)\eta I_{n-4,2}(h)+2\eta^{n-2}(h+\frac{1}{2\eta})^{\frac{n-2}{2}}\big]\\
                  \frac{1}{2n-4}\big[nI_{n-2,1}(h)-n\eta I_{n-2,0}(h)+2\eta^{n-2}(h+\frac{1}{2\eta})^{\frac{n}{2}}\big]
                \end{matrix}\right),
\end{eqnarray}
where
\begin{eqnarray*}
\mathbf{A}=\left(\begin{matrix}
                1&-\frac{1}{h}&0&\cdots&0&0&0\\
                 0&1&0&\cdots&0&0&0\\
                 0&0&1&\cdots&0&0&0\\
                 \vdots&\vdots&\vdots&\vdots&\vdots&\vdots&\vdots\\
                  0&0&0&\cdots&0&1&0\\
                  0&0&0&\cdots&0&0&1
                \end{matrix}\right).\ \
\end{eqnarray*}
Hence, the first equality in (2.5) holds.

 Next we discuss the degrees of the polynomials $\tilde{\alpha}_1(h)$, $\tilde{\beta}_1(h)$ and $\tilde{\varphi}_l(h)$ in (2.5). If $(i,j)=(2,n-2),(4,n-4),\cdots,(n-2,2),(n,0)$, then, in view of (2.14) and noting that $n$ is an even number, we obtain
\begin{eqnarray*}
\begin{aligned}
I_{i,j}(h)=&h\Big[\alpha^{(n-1)}(h)I_{0,1}(h)+\beta^{(n-1)}(h)I_{2,0}(h)+\varphi^{(n-1)}(h)\\
&+\alpha^{(n-2)}(h)I_{0,1}(h)+\beta^{(n-2)}(h)I_{2,0}(h)+\varphi^{(n-2)}(h)+\xi_{\frac{n}{2}}(h)\Big]\\
:=&\alpha^{(n)}(h)I_{0,1}(h)+\beta^{(n)}(h)I_{2,0}(h)+\varphi^{(n)}(h),
\end{aligned}
\end{eqnarray*}
where $\alpha^{(n-s)}(h)$ and $\beta^{(n-s)}(h)\ (s=1,2)$ are polynomials in $h$ satisfying
$$\deg\alpha^{(n-1)}(h)\leq n-2,\ \deg\beta^{(n-1)}(h)\leq n-3,\ \deg\alpha^{(n-2)}(h),\deg\beta^{(n-2)}(h)\leq n-4,$$
 $\varphi^{(n-1)}(h)$  is a polynomial in $h$ satisfying $\deg\varphi^{(n-1)}(h)\leq \frac{3}{2}n-3$, $\varphi^{(n-2)}(h)$  is a polynomial in $h$ satisfying $\deg\varphi^{(n-2)}(h)\leq \frac{3}{2}n-5$ and  $\xi_{\frac{n}{2}}(h)$ is a polynomial in $h$ with degree at most $\frac{n}{2}$. Therefore,
$$\deg\alpha^{(n)}(h), \deg\beta^{(n)}(h)\leq n-2,\ \deg\varphi^{(n)}(h)\leq\frac{3}{2}n-2.$$
Similarly, the conclusion holds for $(i,j)=(0,n)$.

If $n$ is an odd number, we can prove the conclusion similarly. This ends the proof.\quad $\lozenge$
 \vskip 0.2 true cm

From (2.4) and Lemma 2.1, we obtain the algebraic structure of the first Melnikov function $M(h)$ immediately.
\vskip 0.2 true cm

\noindent
{\bf Lemma 2.2.}\, {\it If $i+j=n>3$, then
\begin{eqnarray}
\begin{aligned}
M(h)=&\frac{1}{h^{n-2}}\Big[{\alpha}_1(h)I_{0,1}(h)+{\beta}_1(h)I_{2,0}(h)+{\gamma}_1(h)I_{1,0}(h)+{\delta}_1(h)I_{1,1}(h)\\
&+{\alpha}_2(h)J_{0,1}(h)+{\beta}_2(h)J_{2,0}(h)+{\gamma}_2(h)J_{1,0}(h)+{\delta}_2(h)J_{1,1}(h)\\
&+{\varphi}_{\frac{3}{2}n-\frac{7+(-1)^n}{4}}(h)+\sqrt{ h+\frac{1}{2\eta}}{\psi}_{\frac{3}{2}n-\frac{9-(-1)^n}{4}}(h)\Big],
\end{aligned}
\end{eqnarray}
where ${\varphi}_{l}(h)$ and ${\psi}_{l}(h)$ are polynomials in $h$ of degrees at most $l$, and ${\alpha}_k(h)$, ${\beta}_k(h)$, ${\gamma}_k(h)$ and ${\delta}_k(h)$ are polynomials of $h$  with
\begin{eqnarray*}
\begin{aligned}
&\deg{\alpha}_k(h)\leq n-\frac{3+(-1)^n}{2},\ \deg{\delta}_k(h)\leq n-\frac{3-(-1)^n}{2},\\
&\deg{\beta}_k(h),\deg{\gamma}_k(h)\leq n-2,\ k=1,2.
\end{aligned}\end{eqnarray*}
If $n=1,2,3$, then
\begin{eqnarray}
\begin{aligned}
M(h)=&\frac{1}{h}\Big[{\alpha}_1(h)I_{0,1}(h)+{\beta}_1(h)I_{2,0}(h)+{\gamma}_1(h)I_{1,0}(h)+{\delta}_1(h)I_{1,1}(h)\\
&+{\alpha}_2(h)J_{0,1}(h)+{\beta}_2(h)J_{2,0}(h)+{\gamma}_2(h)J_{1,0}(h)+{\delta}_2(h)J_{1,1}(h)\\
&+{\varphi}_{3}(h)+\sqrt{ h+\frac{1}{2\eta}}{\psi}_{2}(h)\Big],
\end{aligned}
\end{eqnarray}
where ${\varphi}_{l}(h)$ and ${\psi}_{l}(h)$ are polynomials in $h$ of degrees at most $l$, and ${\alpha}_k(h)$, ${\beta}_k(h)$, ${\gamma}_k(h)$ and ${\delta}_k(h)$ are polynomials of $h$  with
\begin{eqnarray*}
\begin{aligned}
\deg{\alpha}_k(h), \deg{\delta}_k(h)\leq 2,\ \
\deg{\beta}_k(h),\deg{\gamma}_k(h)\leq 1,\ k=1,2.
\end{aligned}\end{eqnarray*}}

The following lemma gives the Picard-Fuchs equations which the generators of $M(h)$ satisfy.
\vskip 0.2 true cm
\noindent
{\bf Lemma 2.3.}\, (i) {\it The vector functions $\big(I_{0,1}(h),I_{2,0}(h)\big)^T$ and $\big(I_{1,0}(h),I_{1,1}(h)\big)^T$ respectively satisfy the Picard-Fuchs equations
\begin{eqnarray}
\left(\begin{matrix}
                I_{0,1}(h)\\
                 I_{2,0}(h)\\
                 \end{matrix}\right)
=\left(\begin{matrix}
                2\big(h+\frac{1}{2\eta}\big)&0\\
                 h+\frac{1}{2\eta}&h\\
                  \end{matrix}\right)
                \left(\begin{matrix}
                I'_{0,1}(h)\\
                 I'_{2,0}(h)\\
                 \end{matrix}\right)+
                \left(\begin{matrix}
                0\\
                 -\frac{1}{2}(h+\frac{1}{2\eta})\\
                \end{matrix}\right)
\end{eqnarray}
and
\begin{eqnarray}
\left(\begin{matrix}
                 I_{1,0}(h)\\
                 I_{1,1}(h)\\
          \end{matrix}\right)
=\left(\begin{matrix}
                h+\frac{1}{2\eta}&0\\
                 1&2h\\
                \end{matrix}\right)
                \left(\begin{matrix}
                I'_{1,0}(h)\\
                 I'_{1,1}(h)\\
                \end{matrix}\right)+
                \left(\begin{matrix}
                0\\
                  -\sqrt{h+\frac{1}{2\eta}}\\
                \end{matrix}\right).
\end{eqnarray}}
(ii) {\it The vector functions $\big(J_{0,1}(h),J_{2,0}(h)\big)^T$ and $\big(J_{1,0}(h),J_{1,1}(h)\big)^T$ respectively satisfy the Picard-Fuchs equations
\begin{eqnarray}
\left(\begin{matrix}
                J_{0,1}(h)\\
                 J_{2,0}(h)\\
                 \end{matrix}\right)
=\left(\begin{matrix}
                2\big(h+\frac{1}{2\eta}\big)&0\\
                 h+\frac{1}{2\eta}&h\\
                  \end{matrix}\right)
                \left(\begin{matrix}
                J'_{0,1}(h)\\
                 J'_{2,0}(h)\\
                 \end{matrix}\right)+
                \left(\begin{matrix}
                0\\
                 \frac{1}{2}(h+\frac{1}{2\eta})\\
                \end{matrix}\right)
\end{eqnarray}
and
\begin{eqnarray}
\left(\begin{matrix}
                 J_{1,0}(h)\\
                 J_{1,1}(h)\\
          \end{matrix}\right)
=\left(\begin{matrix}
                h+\frac{1}{2\eta}&0\\
                 1&2h\\
                \end{matrix}\right)
                \left(\begin{matrix}
                J'_{1,0}(h)\\
                 J'_{1,1}(h)\\
                \end{matrix}\right)+
                \left(\begin{matrix}
                0\\
                  \sqrt{h+\frac{1}{2\eta}}\\
                \end{matrix}\right).
\end{eqnarray}}

\noindent
{\bf Proof.}\, We only prove the conclusion (i). Conclusion (ii) can be proved similarly. Since $x$ can be regarded as a function of $y$ and $h$, differentiating the first equation in (1.3) with respect to $h$,  we get $$\frac{\partial x}{\partial h}=\frac{y^2}{2x},$$ which implies
\begin{eqnarray}
I'_{i,j}(h)=\frac{i}{2}\int_{L^1_h}x^{i-2}y^{j-1}dx.
\end{eqnarray}
Hence,
\begin{eqnarray}
I_{i,j}(h)=\frac{2}{i+2}I'_{i+2,j-2}(h).
\end{eqnarray}
Multiplying both side of (2.21) by $h$ and integrating over $L^1_h$, we have
\begin{eqnarray}
hI'_{i,j}(h)=\frac{i}{i+2}I'_{i+2,j-2}(h)-I'_{i,j-1}(h)+\frac{\eta}{2}I'_{i,j-2}(h).
\end{eqnarray}

On the other hand, by (2.2), we have for $i\geq1$ and $j\geq-1$
\begin{eqnarray}
\begin{aligned}
I_{i,j}(h)=&\int_{L^1_h}x^{i}y^{j-3}dy=-\frac{i}{j-2}\int_{L^1_h}x^{i-1}y^{j-2}dx-\frac{i}{j-2}\eta^{j-2}\int_{\overrightarrow{BG}}x^{i-1}dx\\
=&-\frac{i}{2(j-2)}\int_{L^1_h}x^{i-2}y^{j-2}(2hy+1)dy+\frac{\eta^{i+j-2}}{j-2}\big(h+\frac{1}{2\eta}\big)^\frac{i}{2}\\
=&-\frac{1}{j-2}\Big[2hI'_{i,j}(h)+I'_{i,j-1}(h)-\eta^{i+j-2}\big(h+\frac{1}{2\eta}\big)^\frac{i}{2}\Big].
\end{aligned}
\end{eqnarray}
 Taking $(i,j)=(0,1)$ in (2.22) and noting that (2.12) we obtain
 $$I_{0,1}(h)=2\big(h+\frac{1}{2\eta}\big)I'_{0,1}(h).$$
 From (2.24) we have
\begin{eqnarray*}\left\{\begin{array}{l}
I_{2,0}(h)=hI'_{2,0}(h)+\frac{1}{2}I'_{2,-1}(h)-\frac{1}{2}\big(h+\frac{1}{2\eta}\big),\\
I_{1,0}(h)=hI'_{1,0}(h)+\frac{1}{2}I'_{1,-1}(h)-\frac{1}{2\eta}\sqrt{h+\frac{1}{2\eta}},\\
I_{1,1}(h)=2hI'_{1,1}(h)+I'_{1,0}(h)-\sqrt{h+\frac{1}{2\eta}}.
\end{array}\right.
\end{eqnarray*}
In view of (2.9) and (2.12) we obtain the conclusion (i). The proof is completed.\quad $\lozenge$

 \vskip 0.2 true cm

\noindent
{\bf Lemma 2.4.}\, {\it For $h\in\Sigma=(-\frac{1}{2\eta},0)$, we have
\begin{eqnarray}
\begin{cases}
I_{0,1}(h)=-\sqrt{\frac{2}{\eta}}\sqrt{h+\frac{1}{2\eta}},\ \ I_{1,0}(h)=c_1\big(h+\frac{1}{2\eta}\big),\\
I_{2,0}(h)=\frac{1}{2}h\ln\frac{1-\sqrt{ 2\eta h+1}}{1+\sqrt{2\eta h+1}}+\frac{1}{2}h\ln|h|-\frac{1}{2\eta}\sqrt{2\eta h+1}-c_2h-\frac{1}{4\eta},\\
I_{1,1}(h)=\frac{1}{2}\sqrt{|h|}\arctan\frac{2\eta h+\frac{1}{2}}{\sqrt{-2\eta h(2\eta h+1)}}-\sqrt{ h+\frac{1}{ 2\eta}}+\big(\frac{\pi}{4}-c_1\sqrt{2\eta}\big)\sqrt{|h|}+c_1
\end{cases}
\end{eqnarray}
and
\begin{eqnarray}
\begin{cases}
J_{0,1}(h)=-\sqrt{\frac{2}{\eta}}\sqrt{h+\frac{1}{2\eta}},\ \ J_{1,0}(h)=d_1\big(h+\frac{1}{2\eta}\big),\\
J_{2,0}(h)=\frac{1}{2}h\ln\frac{1-\sqrt{2\eta h+1}}{1+\sqrt{2\eta h+1}}-\frac{1}{2}h\ln|h|-\frac{1}{2\eta}\sqrt{2\eta h+1}-d_2h+\frac{1}{4\eta},\\
J_{1,1}(h)=-\frac{1}{2}\sqrt{|h|}\arctan\frac{2\eta h+\frac{1}{2}}{\sqrt{-2\eta h(2\eta h+1)}}+\sqrt{ h+\frac{1}{ 2\eta}}-\big(\frac{\pi}{4}+d_1\sqrt{2\eta}\big)\sqrt{|h|}+d_1,
\end{cases}
\end{eqnarray}
where $c_i$ and $d_i$ ($i=1,2$) are real constants.}
\vskip 0.2 true cm

\noindent
{\bf Proof.} We only prove (2.25). (2.26) can be shown in a similar way. Since the coordinates of $A$ and $C$ are $(0,\frac{-\sqrt{2\eta h+1}-1}{2h})$ and $(0,\frac{\sqrt{2\eta h+1}-1}{2h})$, we have
 \begin{eqnarray*}
I_{0,1}(h)=\int_{L^1_h}y^{-2}dy=-\sqrt{\frac{2}{\eta}}\sqrt{h+\frac{1}{2\eta}}.
\end{eqnarray*}
It follows from the first equation in (2.18) that
\begin{eqnarray*}
I_{1,0}(h)=c_1\big(h+\frac{1}{2\eta}\big),
\end{eqnarray*}
where $c_1$ is a real constant. In the meanwhile, we have
\begin{eqnarray*}
\begin{aligned}
I_{2,0}(h)=\frac{1}{2}h\ln\frac{1-\sqrt{2\eta h+1}}{1+\sqrt{2\eta h+1}}+\frac{1}{2}h\ln|h|-\frac{1}{2\eta}\sqrt{2\eta h+1}-c_2h-\frac{1}{4\eta},\\
I_{1,1}(h)=\frac{1}{2}\sqrt{|h|}\arctan\frac{2\eta h+\frac{1}{2}}{\sqrt{-2\eta h(2\eta h+1)}}-\sqrt{ h+\frac{1}{ 2\eta}}+c\sqrt{|h|}+c_1,
\end{aligned}
\end{eqnarray*}
where $c$ and $c_2$ are real constants. Since $I_{1,1}(\frac{1}{-2\eta})=0$, we have $c=\frac{\pi}{4}-c_1\sqrt{2\eta}$. This completes the proof.\quad $\lozenge$
\section{Proof of the Theorem 1.1}
 \setcounter{equation}{0}
\renewcommand\theequation{3.\arabic{equation}}

In the following, we denote by $P_k(u)$, $Q_k(u)$, $R_k(u)$, $S_k(u)$ and $T_k(u)$ the polynomials of $u$ with degree at most $k$ and denote by $\#\{\phi(h)=0, h\in(\lambda_1,\lambda_2)\}$ the number of isolated zeros of $\phi(h)$ on $(\lambda_1,\lambda_2)$ taking into account the multiplicity.
\vskip 0.2 true cm

\noindent
{\bf Proof of the Theorem 1.1.} If $n>3$ is an even number, let $M_1(h)=h^{n-2}M(h)$ for $h\in(-\frac{1}{2\eta},0)$, then $M_1(h)$ and $M(h)$ have the same number of zeros on $(-\frac{1}{2\eta},0)$.  By Lemmas 2.2 and 2.4, we have
\begin{eqnarray*}
\begin{aligned}
M_1(h)=&{\alpha}_1(h)I_{0,1}(h)+{\beta}_1(h)I_{2,0}(h)+{\gamma}_1(h)I_{1,0}(h)+{\delta}_1(h)I_{1,1}(h)\\
&+{\alpha}_2(h)J_{0,1}(h)+{\beta}_2(h)J_{2,0}(h)+{\gamma}_2(h)J_{1,0}(h)+{\delta}_2(h)J_{1,1}(h)\\
&+{\varphi}_{\frac{3}{2}n-\frac{7+(-1)^n}{4}}(h)+\sqrt{ h+\frac{1}{2\eta}}{\psi}_{\frac{3}{2}n-\frac{9-(-1)^n}{4}}(h)\\
:=&P_{n-1}(h)\ln\frac{1-\sqrt{ 2\eta h+1}}{1+\sqrt{2\eta h+1}}+Q_{n-1}(h)\sqrt{|h|}\arctan\frac{2\eta h+\frac{1}{2}}{\sqrt{-2\eta h(2\eta h+1)}}\\
&+R_{n-1}(h)\ln|h|+S_{n-1}(h)\sqrt{|h|}+{\varphi}_{\frac{3}{2}n-2}(h)+\sqrt{ h+\frac{1}{2\eta}}{\psi}_{\frac{3}{2}n-2}(h).
\end{aligned}
\end{eqnarray*}
Let $t=\sqrt{h+\frac{1}{2\eta}}$, $t\in (0,\frac{1}{\sqrt{2\eta}})$, we have
$$\begin{aligned}
M_1(t)=&P_{n-1}(t^2)\ln\frac{1-\sqrt{ 2\eta}t}{1+\sqrt{2\eta}t}+Q_{n-1}(t^2)\sqrt{\frac{1}{2\eta}-t^2}\arctan\frac{2\eta t^2-\frac{1}{2}}{\sqrt{2\eta (1-2\eta t^2)}t}\\
&+R_{n-1}(t^2)\ln\big(\frac{1}{2\eta}-t^2\big)+S_{n-1}(t^2)\sqrt{\frac{1}{2\eta}-t^2}+T_{3n-3}(t).\end{aligned}$$
Hence, $M_1(h)$ and $M_1(t)$ have the same number of zeros for $h\in(-\frac{1}{2\eta})$ and $t\in(0,\frac{1}{\sqrt{2\eta}})$. Suppose that $\Sigma_1=(0,\frac{1}{\sqrt{2\eta}})\backslash \{t\in(0,\frac{1}{\sqrt{2\eta}})|P_{n-1}(t^2)=0\}$. Then, for $h\in \Sigma_1$, we get
{\footnotesize\begin{eqnarray*}
\begin{aligned}
&\frac{d}{dt}\Big(\frac{M_1(t)}{P_{n-1}(t^2)}\Big)=\frac{2\sqrt{2\eta}}{2\eta t^2-1}+\\
&\frac{d}{dt}\Big[\frac{Q_{n-1}(t^2)\sqrt{\frac{1}{2\eta}-t^2}\arctan\frac{2\eta t^2-\frac{1}{2}}{\sqrt{2\eta (1-2\eta t^2)}t}
+R_{n-1}(t^2)\ln\big(\frac{1}{2\eta}-t^2\big)+S_{n-1}(t^2)\sqrt{\frac{1}{2\eta}-t^2}+T_{3n-3}(t)}{P_{n-1}(t^2)}\Big]\\
=&\frac{P_{2n-2}(t^2)t\ln\big(\frac{1}{2\eta}-t^2\big)+Q_{2n-2}(t^2)t\sqrt{1-2\eta t^2}\arctan\frac{2\eta t^2-\frac{1}{2}}{\sqrt{2\eta (1-2\eta t^2)}t}  +R_{2n-2}(t^2)t\sqrt{1-2\eta t^2}+S_{5n-4}(t)}{P^2_{n-1}(t^2)(2\eta t^2-1)}\\
:=&\frac{M_2(t)}{P^2_{n-1}(t^2)(2\eta t^2-1)}.
\end{aligned}
\end{eqnarray*}}
Let $\Sigma_2=(0,\frac{1}{\sqrt{2\eta}})\backslash \{t\in(0,\frac{1}{\sqrt{2\eta}})|P_{2n-2}(t^2)=0\}$, then we have for $h\in\Sigma_2$
 {\footnotesize\begin{eqnarray*}
\begin{aligned}
\frac{d}{dt}\Big(\frac{M_2(t)}{P_{2n-2}(t^2)t}\Big)=&\frac{d}{dt}\Big[\frac{Q_{2n-2}(t^2)t\sqrt{1-2\eta t^2}\arctan\frac{2\eta t^2-\frac{1}{2}}{\sqrt{2\eta (1-2\eta t^2)}t}  +R_{2n-2}(t^2)t\sqrt{1-2\eta t^2}+S_{5n-4}(t)}{P_{2n-2}(t^2)t}\Big]\\
&+\frac{4\eta t}{2\eta t^2-1}\\
=&\frac{P_{4n-3}(t^2)t\sqrt{1-2\eta t^2}\arctan\frac{2\eta t^2-\frac{1}{2}}{\sqrt{2\eta (1-2\eta t^2)}t}+Q_{4n-3}(t^2)t\sqrt{1-2\eta t^2}+R_{9n-6}(t)}{P^2_{2n-2}(t^2)t^2(2\eta t^2-1)}\\
:=&\frac{M_3(t)}{P^2_{2n-2}(t^2)t^2(2\eta t^2-1)}.
\end{aligned}
\end{eqnarray*}}
Similarly, let $\Sigma_3=(0,\frac{1}{\sqrt{2\eta}})\backslash \{t\in(0,\frac{1}{\sqrt{2\eta}})|P_{4n-3}(t^2)=0\}$, then we have for $h\in\Sigma_3$
 {\footnotesize\begin{eqnarray*}
\begin{aligned}
\frac{d}{dt}\Big(\frac{M_3(t)}{P_{4n-3}(t^2)t\sqrt{1-2\eta t^2}}\Big)=&\frac{d}{dt}\Big[\frac{Q_{4n-3}(t^2)t\sqrt{1-2\eta t^2}+R_{9n-6}(t)}{P_{4n-3}(t^2)t\sqrt{1-2\eta t^2}}\Big]+\frac{2\sqrt{2\eta t}}{\sqrt{1-2\eta t^2}}\\
=&\frac{P_{17n-10}(t)+Q_{8n-5}(t^2)t\sqrt{1-2\eta t^2}}{P^2_{4n-3}(t^2)t^2(1-2\eta t^2)\sqrt{1-2\eta t^2}}\\
:=&\frac{M_4(t)}{P^2_{4n-3}(t^2)t^2(1-2\eta t^2)\sqrt{1-2\eta t^2}}.
\end{aligned}
\end{eqnarray*}}
Let $M_4(t)=P_{17n-10}(t)+Q_{8n-5}(t^2)t\sqrt{1-2\eta t^2}=0$. That is,
$$Q_{8n-5}(t^2)t\sqrt{1-2\eta t^2}=-P_{17n-10}(t).$$
By squaring the above equation, we can deduce that $M_4(t)$ has at most $34n-20$ zeros on $(0,\frac{1}{\sqrt{2\eta}})$.
Hence,
 \begin{eqnarray}
 \begin{aligned}
 \# \{M(h)=0, h\in (-\frac{1}{2\eta},0)\}&=\# \{M_1(t)=0, t\in (0,\frac{1}{\sqrt{2\eta}})\}\\
&\le 41n-23.
\end{aligned}
\end{eqnarray}
If $n=1,2,3$, it is easy to check that (3.1) also holds.

If $n$ is an odd number, we can prove the Theorem 1.1 similarly. This ends the proof of Theorem 1.1. \quad $\lozenge$

\section{Proof of the Theorem 1.2}
 \setcounter{equation}{0}
\renewcommand\theequation{4.\arabic{equation}}

If $f^1(x,y)=f^2(x,y)$ and $f^3(x,y)=f^4(x,y)$, then system (1.2) can be written as
\begin{eqnarray}
\left(  \begin{array}{c}
          \dot{x} \\          \dot{y}
          \end{array} \right)
=\begin{cases}
 \left(
  \begin{array}{c}
          y-2x^2-\eta+\varepsilon f^1(x,y) \\
          -2xy+\varepsilon g^1(x,y)
          \end{array} \right), \quad x>0,\\[0.6truecm]
  \left(  \begin{array}{c}
         y-2x^2-\eta+\varepsilon f^3(x,y) \\
          -2xy+\varepsilon g^3(x,y)
           \end{array}
 \right),\quad x<0.
  \end{cases}
 \end{eqnarray}
From Theorem 1.1 in \cite{LH}, we know that the first order Melnikov function $M(h)$ of system (4.1) has the following form
\begin{eqnarray}
\begin{aligned}
M(h)=&\frac{H^1_y(A)H^3_y(C)}{H^3_y(A)H^1_y(C)}\int_{L^1_h}\mu^1(x,y)[g^1(x,y)dx-f^1(x,y)dy]\\
&+\frac{H^1_y(A)}{H^3_y(A)}
\int_{L^3_h}\mu^3(x,y)[g^3(x,y)dx-f^3(x,y)dy],\ h\in\Sigma
\end{aligned}
\end{eqnarray}
and the number of zeros of $M(h)$ controls the number of limit cycles of system (4.1) if $M(h)\not\equiv0$
in the corresponding period annulus. Noting that
$\frac{H^1_y(A)H^3_y(C)}{H^3_y(A)H^1_y(C)}=\frac{H^1_y(A)}{H^3_y(A)}=1,$ we have
\begin{eqnarray*}
\begin{aligned}
M(h)=\int_{L^1_h}y^{-3}[g^1(x,y)dx-f^1(x,y)dy]+\int_{L^3_h}y^{-3}[g^3(x,y)dx-f^3(x,y)dy].
\end{aligned}
\end{eqnarray*}

For $h\in\Sigma$ and $i=0,1,2,\cdots,j=-1,0,1,2,\cdots$, we denote
$$\begin{aligned}
U_{i,j}(h)=\int_{\Gamma_h}x^iy^{j-3}dy,\ \ \tilde{U}_{i,j}(h)=\int_{\tilde{\Gamma}_h}x^iy^{j-3}dy,
\end{aligned}$$
where $\Gamma_h=L^1_h\cup L^2_h$ and $\tilde{\Gamma}_h=L^3_h\cup L^4_h$. It is easy to get that $\tilde{U}_{i,j}(h)=(-1)^{i+1}{U}_{i,j}(h)$. Similar to (2.2), we get
$$\begin{aligned}&\int_{\Gamma_h}x^iy^jdx=-\frac{j}{i+1}\int_{\Gamma_h}x^{i+1}y^{j-1}dy,\\ &\int_{\tilde{\Gamma}_h}x^iy^jdx=-\frac{j}{i+1}\int_{\tilde{\Gamma}_h}x^{i+1}y^{j-1}dy.\end{aligned}$$
Therefore,
\begin{eqnarray*}
\begin{aligned}
M(h)=&\sum\limits_{i+j=0}^nb^1_{i,j}\int_{\Gamma_h}x^iy^{j-3}dx
-\sum\limits_{i+j=0}^na^1_{i,j}\int_{\Gamma_h}x^iy^{j-3}dy\\
&+\sum\limits_{i+j=0}^nb^3_{i,j}\int_{\tilde{\Gamma}_h}x^iy^{j-3}dx
-\sum\limits_{i+j=0}^na^3_{i,j}\int_{\tilde{\Gamma}_h}x^iy^{j-3}dy\\
=&-\sum\limits_{i+j=0}^nb^1_{i,j}\frac{j-3}{i+1}\int_{\Gamma_h}x^{i+1}y^{j-4}dy
-\sum\limits_{i+j=0}^na^1_{i,j}\int_{\Gamma_h}x^iy^{j-3}dy\\
&-\sum\limits_{i+j=0}^nb^3_{i,j}\frac{j-3}{i+1}\int_{\tilde{\Gamma}_h}x^{i+1}y^{j-4}dy
-\sum\limits_{i+j=0}^na^3_{i,j}\int_{\tilde{\Gamma}_h}x^iy^{j-3}dy\\
:=&\sum\limits_{\substack{i+j=0,\\i\geq0,j\geq-1}}^n\bar{\sigma}_{i,j}U_{i,j}(h),
\end{aligned}
\end{eqnarray*}
where $\bar{\sigma}_{i,j}$ are real constants.
\vskip 0.2 true cm
Following the proof of Lemmas 2.1 and 2.2, we can get the following Lemma 4.1.
\vskip 0.2 true cm

\noindent
{\bf Lemma 4.1.}\, {\it If $i+j=n>3$, then
\begin{eqnarray}
\begin{aligned}
M(h)=\frac{1}{h^{n-2}}\Big[\alpha_1(h)U_{0,1}(h)+{\beta_1}(h)U_{2,0}(h)+{\gamma_1}(h)U_{1,0}(h)+{\delta_1}(h)U_{1,1}(h)\Big],
\end{aligned}
\end{eqnarray}
where ${\alpha_1}(h)$, ${\beta_1}(h)$, ${\gamma_1}(h)$ and ${\delta_1}(h)$ are polynomials of $h$  with
\begin{eqnarray*}
\begin{aligned}
&\deg{\alpha_1}(h)\leq n-\frac{3+(-1)^n}{2},\ \deg{\delta_1}(h)\leq n-\frac{3-(-1)^n}{2},\\
&\deg{\beta_1}(h),\deg{\gamma_1}(h)\leq n-2.
\end{aligned}\end{eqnarray*}
If $n=1,2,3$, then
\begin{eqnarray}
\begin{aligned}
M(h)=\frac{1}{h}\Big[{\alpha_1}(h)U_{0,1}(h)+{\beta_1}(h)U_{2,0}(h)+{\gamma_1}(h)U_{1,0}(h)+{\delta_1}(h)U_{1,1}(h)\Big],
\end{aligned}
\end{eqnarray}
where ${\alpha_1}(h)$, ${\beta_1}(h)$, ${\gamma_1}(h)$ and ${\delta_1}(h)$ are polynomials of $h$  with
\begin{eqnarray*}
\begin{aligned}
\deg{\alpha}_1(h), \deg{\delta}_1(h)\leq 2,\ \
\deg{\beta}_1(h),\deg{\gamma}_1(h)\leq 1.
\end{aligned}\end{eqnarray*}}

The following lemma gives the Picard-Fuchs equations which the generators of $M(h)$ in (4.2) satisfy and can be proved by the method in Lemma 2.3.
\vskip 0.2 true cm
\noindent
{\bf Lemma 4.2.}\,{\it The vector functions $\big(U_{0,1}(h),U_{2,0}(h)\big)^T$ and $\big(U_{1,0}(h),U_{1,1}(h)\big)^T$ respectively satisfy the Picard-Fuchs equations
\begin{eqnarray}
\left(\begin{matrix}
                U_{0,1}(h)\\
                U_{2,0}(h)\\
                 \end{matrix}\right)
=\left(\begin{matrix}
                2\big(h+\frac{1}{2\eta}\big)&0\\
                 h+\frac{1}{2\eta}&h\\
                  \end{matrix}\right)
                \left(\begin{matrix}
                U'_{0,1}(h)\\
                 U'_{2,0}(h)\\
                 \end{matrix}\right)
\end{eqnarray}
and
\begin{eqnarray}
\left(\begin{matrix}
                 U_{1,0}(h)\\
                 U_{1,1}(h)\\
          \end{matrix}\right)
=\left(\begin{matrix}
                h+\frac{1}{2\eta}&0\\
                 1&2h\\
                \end{matrix}\right)
                \left(\begin{matrix}
                U'_{1,0}(h)\\
                 U'_{1,1}(h)\\
                \end{matrix}\right).
\end{eqnarray}}

From (4.5) and (4.6), we have for $h\in\Sigma=(-\frac{1}{2\eta},0)$
\begin{eqnarray}
\begin{cases}
U_{0,1}(h)=-{\frac{2}{\eta}}\sqrt{h+\frac{1}{2\eta}},\ \ U_{1,0}(h)=e_1\big(h+\frac{1}{2\eta}\big),\\
U_{2,0}(h)=\frac{1}{2}h\ln\frac{1-\sqrt{ 2\eta h+1}}{1+\sqrt{2\eta h+1}}-\frac{1}{2\eta}\sqrt{2\eta h+1}-e_2h,\\
U_{1,1}(h)=\big(\frac{\pi}{4}-e_1\sqrt{2\eta}\big)\sqrt{|h|}+e_1,
\end{cases}
\end{eqnarray}
where $e_1$ and $e_2$ are real constants. Hence,
$$\begin{aligned}
M(h)=&\frac{1}{h^{n-2}}\big[{\alpha_1}(h)U_{0,1}(h)+{\beta_1}(h)U_{2,0}(h)+{\gamma_1}(h)U_{1,0}(h)+{\delta_1}(h)U_{1,1}(h)\big]\\
=&\frac{1}{h^{n-2}}\Big[-{\frac{2}{\eta}}\alpha_1(h)\sqrt{h+\frac{1}{2\eta}}+\beta_1(h)\Big(\frac{1}{2}h\ln\frac{1-\sqrt{2\eta h+1}}{1+\sqrt{2\eta h+1}}-\frac{1}{2\eta}\sqrt{2\eta h+1}-e_2h\Big)\\
&+e_1\gamma_1(h)\big(h+\frac{1}{2\eta}\big)+\delta_1(h)\big((\frac{\pi}{4}-e_1\sqrt{2\eta})\sqrt{|h|}+e_1\big)\Big]\\
:=&P_{n-1}(h)\ln\frac{1-\sqrt{2\eta h+1}}{1+\sqrt{2\eta h+1}}+Q_{n-2}(h)\sqrt{2\eta h+1}+R_{n-1}(h)\sqrt{|h|}+S_{n-1}(h).
\end{aligned}$$
Following the lines of the proof of Theorem 1.1, we can prove that $M(h)$ has at most $9n-4$ zeros on $(-\frac{1}{2\eta},0)$. The Theorem 1.2 is proved.

\section{Proof of the Theorem 1.3}
 \setcounter{equation}{0}
\renewcommand\theequation{5.\arabic{equation}}

If $f^1(x,y)=f^4(x,y)$ and $f^2(x,y)=f^3(x,y)$, then system (1.2) can be written as
\begin{eqnarray}
\left(  \begin{array}{c}
          \dot{x} \\          \dot{y}
          \end{array} \right)
=\begin{cases}
 \left(
  \begin{array}{c}
          y-2x^2-\eta+\varepsilon f^1(x,y) \\
          -2xy+\varepsilon g^1(x,y)
          \end{array} \right), \quad y>\eta,\\[0.6truecm]
  \left(  \begin{array}{c}
         y-2x^2-\eta+\varepsilon f^2(x,y) \\
          -2xy+\varepsilon g^2(x,y)
           \end{array}
 \right),\quad y<\eta.
  \end{cases}
 \end{eqnarray}
From Theorem 1.1 in \cite{LH}, we know that the first order Melnikov function $M(h)$ of system (5.1) has the following form
\begin{eqnarray}
\begin{aligned}
M(h)=&\frac{H^1_x(D)H^2_x(B)}{H^2_x(D)H^1_x(B)}\int_{\Upsilon_h}\mu^1(x,y)[g^1(x,y)dx-f^1(x,y)dy]\\
&+\frac{H^1_x(D)}{H^2_x(D)}
\int_{\tilde{\Upsilon}_h}\mu^2(x,y)[g^2(x,y)dx-f^2(x,y)dy],\ h\in\Sigma
\end{aligned}
\end{eqnarray}
and the number of zeros of $M(h)$ controls the number of limit cycles of system (5.1) if $M(h)\not\equiv0$
in the corresponding period annulus. Noting that
$\frac{H^1_x(D)H^2_x(B)}{H^2_x(D)H^1_x(B)}=\frac{H^1_x(D)}{H^2_x(D)}=1,$ we have
\begin{eqnarray*}
\begin{aligned}
M(h)=\int_{\Upsilon_h}y^{-3}[g^1(x,y)dx-f^1(x,y)dy]+\int_{\tilde{\Upsilon}_h}y^{-3}[g^2(x,y)dx-f^2(x,y)dy].
\end{aligned}
\end{eqnarray*}

For $h\in\Sigma$ and $i=0,1,2,\cdots,j=-1,0,1,2,\cdots$, we denote
$$\begin{aligned}
V_{i,j}(h)=\int_{\Upsilon_h}x^iy^{j-3}dy,\ \ \tilde{V}_{i,j}(h)=\int_{\tilde{\Upsilon}_h}x^iy^{j-3}dy,
\end{aligned}$$
where $\Upsilon_h=L^1_h\cup L^4_h$ and $\tilde{\Upsilon}_h=L^2_h\cup L^3_h$. Noting that $\Upsilon_h$ and $\tilde{\Upsilon}_h$ are symmetric with respect to $x=0$, we get ${V}_{2l,j}(h)=\tilde{V}_{2l,j}(h)=0$ for $l=0,1,2,\cdots$. Similar to (2.2), we have
$$\begin{aligned}&\int_{\Upsilon_h}x^iy^jdx=-\frac{j}{i+1}\int_{\Upsilon_h}x^{i+1}y^{j-1}dy-\eta^j\int_{\overrightarrow{BD}}x^idx,\\ &\int_{\tilde{\Upsilon}_h}x^iy^jdx=-\frac{j}{i+1}\int_{\tilde{\Upsilon}_h}x^{i+1}y^{j-1}dy-\eta^j\int_{\overrightarrow{DB}}x^idx.\end{aligned}$$
Therefore,
\begin{eqnarray*}
\begin{aligned}
M(h)=&\sum\limits_{i+j=0}^nb^1_{i,j}\int_{\Upsilon_h}x^iy^{j-3}dx
-\sum\limits_{i+j=0}^na^1_{i,j}\int_{\Upsilon_h}x^iy^{j-3}dy\\
&+\sum\limits_{i+j=0}^nb^2_{i,j}\int_{\tilde{\Upsilon}_h}x^iy^{j-3}dx
-\sum\limits_{i+j=0}^na^2_{i,j}\int_{\tilde{\Upsilon}_h}x^iy^{j-3}dy\\
=&\sum\limits_{\substack{i+j=0,\\i\geq0,j\geq-1}}^n\tilde{\sigma}_{i,j}V_{i,j}(h)+
\sum\limits_{\substack{i+j=0,\\i\geq0,j\geq-1}}^n\tilde{\tau}_{i,j}\tilde{V}_{i,j}(h)+
\sum\limits_{i=0}^n\tilde{\nu}_i[(-1)^{i+1}-1]\eta^{i+1}\big(h+\frac{1}{2\eta}\big)^{\frac{i+1}{2}},
\end{aligned}
\end{eqnarray*}
where $\tilde{\sigma}_{i,j}$, $\tilde{\tau}_{i,j}$ and $\tilde{\nu}_i$ are real constants.
\vskip 0.2 true cm

\noindent
{\bf Lemma 5.1.}\, {\it If $i+j=n>3$, then
\begin{eqnarray}
\begin{aligned}
M(h)=\frac{1}{h^{n-2}}\Big[\gamma_1(h)V_{1,0}(h)+{\delta_1}(h)V_{1,1}(h)+\gamma_2(h)\tilde{V}_{1,0}(h)+{\delta_2}(h)\tilde{V}_{1,1}(h)\Big],
\end{aligned}
\end{eqnarray}
where ${\gamma_k}(h)$ and ${\delta_k}(h)$ $(k=1,2)$ are polynomials of $h$  with
\begin{eqnarray*}
\begin{aligned}
\deg{\gamma_k}(h)\leq n-2,\ \ \deg{\delta_k}(h)\leq n-\frac{3-(-1)^n}{2},\ k=1,2.
\end{aligned}\end{eqnarray*}
If $n=1,2,3$, then
\begin{eqnarray}
\begin{aligned}
M(h)=\frac{1}{h}\Big[\gamma_1(h)V_{1,0}(h)+{\delta_1}(h)V_{1,1}(h)+\gamma_2(h)\tilde{V}_{1,0}(h)+{\delta_2}(h)\tilde{V}_{1,1}(h)\Big],
\end{aligned}
\end{eqnarray}
where ${\gamma_k}(h)$ and ${\delta_k}(h)$ $(k=1,2)$ are polynomials of $h$  with
\begin{eqnarray*}
\begin{aligned}
\deg{\gamma}_k(h)\leq1, \ \ \deg{\delta}_k(h)\leq 2,\ k=1,2.
\end{aligned}\end{eqnarray*}}


\noindent
{\bf Lemma 5.2.}\,{\it The vector functions $\big(V_{1,0}(h),V_{1,1}(h)\big)^T$ and $\big(\tilde{V}_{1,0}(h),\tilde{V}_{1,1}(h)\big)^T$ respectively satisfy the Picard-Fuchs equations
\begin{eqnarray}
\left(\begin{matrix}
                V_{1,0}(h)\\
                V_{1,1}(h)\\
                 \end{matrix}\right)
=\left(\begin{matrix}
                h+\frac{1}{2\eta}&0\\
                1&2h\\
                  \end{matrix}\right)
                \left(\begin{matrix}
                V'_{1,0}(h)\\
                 V'_{1,1}(h)\\
                 \end{matrix}\right)+ \left(\begin{matrix}
                0\\
                 -2\sqrt{h+\frac{1}{2\eta}}\\
                 \end{matrix}\right)
\end{eqnarray}
and
\begin{eqnarray}
\left(\begin{matrix}
                \tilde{V}_{1,0}(h)\\
                \tilde{V}_{1,1}(h)\\
                 \end{matrix}\right)
=\left(\begin{matrix}
               h+\frac{1}{2\eta}&0\\
                 1&2h\\
                  \end{matrix}\right)
                \left(\begin{matrix}
                \tilde{V}'_{1,0}(h)\\
                 \tilde{V}'_{1,1}(h)\\
                 \end{matrix}\right)+ \left(\begin{matrix}
                0\\
                 2\sqrt{h+\frac{1}{2\eta}}\\
                 \end{matrix}\right).
\end{eqnarray}}

From (5.5) and (5.6), we have for $h\in\Sigma=(-\frac{1}{2\eta},0)$
\begin{eqnarray}
\begin{aligned}
&V_{1,0}(h)=\hat{c}_1\big(h+\frac{1}{2\eta}\big),\ \ \tilde{V}_{1,0}(h)=\hat{d}_1\big(h+\frac{1}{2\eta}\big),\\
&V_{1,1}(h)=\sqrt{|h|}\arctan\frac{2\eta h+\frac{1}{2}}{\sqrt{-2\eta h(2\eta h+1)}}-2\sqrt{ h+\frac{1}{ 2\eta}}+\big(\frac{\pi}{2}-\hat{c}_1\sqrt{2\eta}\big)\sqrt{|h|}+\hat{c}_1,\\
&\tilde{V}_{1,1}(h)=-\sqrt{|h|}\arctan\frac{2\eta h-\frac{1}{2}}{\sqrt{-2\eta h(2\eta h+1)}}-2\sqrt{ h+\frac{1}{ 2\eta}}-\big(\frac{\pi}{2}+\hat{d}_1\sqrt{2\eta}\big)\sqrt{|h|}+\hat{d}_1,
\end{aligned}
\end{eqnarray}
where $\hat{c}_1$ and $\hat{d}_2$ are real constants. Hence,
{\small$$\begin{aligned}
M(h)=&\frac{1}{h^{n-2}}\big[{\gamma_1}(h)V_{1,0}(h)+{\delta_1}(h)V_{1,1}(h)+{\gamma_2}(h)\tilde{V}_{1,0}(h)+{\delta_2}(h)\tilde{V}_{1,1}(h)\big]\\
=&\frac{1}{h^{n-2}}\Big[{\gamma_1}(h)\hat{c}_1\big(h+\frac{1}{2\eta}\big)+{\gamma_2}(h)\hat{d}_1\big(h+\frac{1}{2\eta}\big)\\
&+{\delta_1}(h)\Big(\sqrt{|h|}\arctan\frac{2\eta h+\frac{1}{2}}{\sqrt{-2\eta h(2\eta h+1)}}-2\sqrt{ h+\frac{1}{ 2\eta}}+\big(\frac{\pi}{2}-\hat{c}_1\sqrt{2\eta}\big)\sqrt{|h|}+\hat{c}_1\Big)\\
&+{\delta_2}(h)\Big(-\sqrt{|h|}\arctan\frac{2\eta h-\frac{1}{2}}{\sqrt{-2\eta h(2\eta h+1)}}-2\sqrt{ h+\frac{1}{ 2\eta}}-\big(\frac{\pi}{2}+\hat{d}_1\sqrt{2\eta}\big)\sqrt{|h|}+\hat{d}_1\Big)\Big]\\
:=&P_{n-1}(h)\sqrt{|h|}\arctan\frac{2\eta h-\frac{1}{2}}{\sqrt{-2\eta h(2\eta h+1)}}+Q_{n-1}(h)\sqrt{2\eta h+1}+R_{n-1}(h)\sqrt{|h|}+S_{n-1}(h).
\end{aligned}$$}
Following the lines of the proof of Theorem 1.1, we get that $M(h)$ has at most $9n-6$ zeros on $(-\frac{1}{2\eta},0)$. The Theorem 1.3 is proved.
\vskip 0.2 true cm

\noindent
{\bf Remark 5.1.} If the switching line is parallel to $y$-axis, then the first order Melnikov function can be written as (4.2). If the switching line is parallel to $x$-axis, then the first order Melnikov function can be written as (5.2).
\vskip 0.2 true cm

\noindent
{\bf Acknowledgment}
 \vskip 0.2 true cm

\noindent
Supported by National Natural Science Foundation of China (11701306,11671040,11601250), Construction of First-class Disciplines of Higher Education of Ningxia(pedagogy)\\ (NXYLXK2017B11), Higher Educational Science Program of Ningxia(NGY201789) and Key Program of Ningxia Normal University(NXSFZD1708).

\end{document}